\def\etal{\mbox{et al.}}
\documentclass[]{gNST2e}
\usepackage{multirow}
\usepackage[pdftex]{hyperref}
\newcommand{\esp}[1]{\mathbb{E} \left( #1 \right)} 
\newcommand{\esploi}[2]{\mathbb{E}_{#2} \left( #1 \right)} 
\newcommand{\var}[1]{\mathbb{V} \left( #1 \right)} 
\newcommand{\varloi}[2]{\mathbb{V}_{#2} \left( #1 \right)} 
\newcommand{\proba}[1]{\mathbb{P} \left( #1 \right)} 
\newcommand{\cov}[1]{\mathbb{C}ov \left( #1 \right)} 
\newcommand{\ind}[1]{\mathbf{1}_{\left\{ #1 \right\}}}
\newcommand{\Ind}[1]{\mathbf{1}{\left\{ #1 \right\}}}
\newcommand{\card}[1]{\left| #1 \right|}
\newcommand{\ent}[1]{\left\lfloor #1 \right\rfloor}
\newcommand{\ens}[1]{\{1,\cdots,#1\}}
\newcommand{\nbt}[2]{#1 \cdot 10^{#2}} 
\newcommand{\ordinal}[1]{$#1^{\textrm{th}}$}
\newcommand{\important}[1]{\textbf{#1}}
\newcommand{\gui}[1]{"{#1}"}
\newcommand{\cod}[1]{\texttt{#1}}

\begin{document}

\doi{10.1080/1048525YYxxxxxxxx}
 \issn{1029-0311}
\issnp{1048-5252}
\jvol{00} \jnum{00} \jyear{2010} 

\markboth{J\'er\^ome Collet}{Journal of Nonparametric Statistics}


\title{Estimate dependence in medium dimensions, using ranks and sub-sampling}

\author{J\'er\^ome Collet$^{\rm a}$$^{\ast}$\thanks{$^\ast$Corresponding author. Email: Jerome.Collet@laposte.net
} \\\vspace{6pt} $^{\rm a}${\em{MODAL'X, Universit\'e Paris Ouest, B\^atiment G, 200 avenue de la R\'epublique, 92000 Nanterre, France}}\\\vspace{6pt}\received{v3.7 released September 2009}}

\maketitle

\begin{abstract}
It is well known that non-parametric methods suffer from the \gui{curse of dimensionality}. We propose here a new estimation method for a multivariate distribution, using sub-sampling and ranks, which seems not to suffer from this \gui{curse}. We prove that in case of independence, the uncertainty of the estimated distribution increases almost linearly w.r.t. the dimension, for dimensions around~6. Otherwise, a simulation study shows that if we use this estimation to build an independence test, the number of observations needed to obtain a given power increases linearly with the dimension. Finally, we give examples of a regression using this estimation: with 3000~observations, in dimension~5, with a markedly complicated dependence, the estimated distribution is graphically very similar to the real one.
\begin{keywords}U-statistics, sub-sampling, ranks, regression\end{keywords}
\begin{classcode}62G07, 62G09, 62G08, 62G30 \end{classcode}
\end{abstract}

\section{Introduction}
The goal of this study is to build non-parametric regression models. It is well known that nowadays, non-parametric statistical methods (if we except machine-learning methods such as random forests) do not allow studying data with dimension larger than~2 or~3. The method we propose here seems to escape the \gui{curse of  dimensionality}.\\
The main object in this paper is the rank space. If we have $m$~observations of a $d$-dimensional random variable, for each observation we replace the value of each component by its rank. We obtain $m$~vectors of the discrete space~$\ens{m}^d$ (of which the total cardinality is~$m^d$). As this space is sparsely populated, it is usual to \gui{smooth} the points of this space. One may for example count the number of points in the intersection of half-spaces~\citep{hoeffding,BKR,deheuvels2}, or more recently~\citep{GR,GQR,Kojadinovic20091137,Beran20071805,Bilodeau2005345}. It is also possible to use kernel methods~\citep{jdf,sca,charp}. We propose in subsection~\ref{testtest} a more specialized state of the art for independence testing, because we use our estimation technique for this purpose, in a simulation study. In any case, these methods are built for continuous data, whereas the vector of ranks is essentially discrete.\\
That is why we try keep the rank space as discrete as possible. With this intention, we propose to use sub-sampling. This type of method is now widely used and studied, we only cite the pioneering work~\citep{breiman1996bagging} and one of the most recent surveys~\citep{buhlmann2012bagging}. These methods are applied to \important{estimators}, so the theoretical results in these papers do not fit exactly the use we propose. Nevertheless, the main conclusion is that sub-sampling related methods are smoothing methods, useful when using discontinuous, non-linear functions.\\
One can assert that sub-sampling, combined with a very discontinuous and non-linear function such as ranking, makes an interesting estimator: even on a small sample, the ranks of the sub-samples densely fill the available space.\\
To measure the accuracy of the estimation, we use the $L^2$~distance between the estimated probability and the theoretical one. One obtains the following results:
\begin{itemize}
  \item In case of independence, for dimensions smaller around~6, the variance of this distance increases almost linearly with the dimension.
  \item Simulations confirm this theoretical result.
  \item If we use these estimations to test independence, the number of observations needed to obtain a given power increases linearly with the dimension (if smaller than~6).
  \item We use this estimation to build a regression model: the dependence studied is non-monotonic, we have 3000~observations, and the dimension is~5. One can see on graphs that the estimation is good.
\end{itemize}

\section{Description of the estimation method}
\subsection{Notation and principles}
The first step is to define the random variable~$\mathbf{X}=(X_1, \cdots, X_d)$ in~$\mathbb{R}^d$, with distribution~$c(\mathbf{u})=c(u_1,\cdots,u_d)$; its marginals are assumed to be continuous.\\
We use an $m$-sample~$\mathbf{X}(\ens{m})$, with $X_l(i)$ the \ordinal{l}~component of the \ordinal{i}~observation of the sample. We note~$\mathbf{R}\left(i,\mathbf{X}(\ens{m})\right)$ the $d$-uple of the ranks of the components of the \ordinal{i}~observation of sample~$\mathbf{X}(\ens{m})$.\\
This allows defining a random variable~$\mathbf{B}$: it is an array, filled with~0s and~1s, $d$~dimensional, each dimension being indexed from~1 to~$m$. For a $d$-uple of ranks~$\mathbf{r} \in \ens{m}^d$, $\mathbf{B}_{\mathbf{r}}(X(\ens{m}))$ is equal to~1 if and only if it is reached by an observation of the $m$-sample. In other words:
\begin{gather*}
  \mathbf{B}_{\mathbf{r}}(X(\ens{m}))=\Ind
  {\exists i \in \ens{m}\ | \ \mathbf{R}\left(i,\mathbf{X}(\ens{m})\right)=\mathbf{r}}
.\end{gather*}
The object we want to estimate is
\begin{gather*}
  P(\mathbf{r},c)=\frac{1}{m}\esploi{\mathbf{B}_{\mathbf{r}}(X(\ens{m}))}{c}
.\end{gather*}
Dividing by~$m$ makes the sum of~$P$ equal to~1, as~$\forall \mathbf{r}, \forall X,  \sum_{\mathbf{r}}\mathbf{B}(\mathbf{r},X(\ens{m}))=m$.\\
We will estimate~$P(\mathbf{r},c)$ from a $n$-sample ($n>m$), using the following $U$-statistic:
\begin{gather*}
  \hat{P}_n(\mathbf{r}) = \frac{1}{m\binom{n}{m}} \sum_{S \in \mathcal{S}}{\mathbf{B}_{\mathbf{r}}(X(S(\ens{m})))}
,\end{gather*}
where~$\mathcal{S}$ is the set of the injections from~$\ens{m}$ to~$\ens{n}$.\\
We will make some remarks:
\begin{itemize}
  \item The estimator~$\hat{P}_n$ can be viewed as a generalization of Kendall's~$\tau$.
  \item A simple example, in 2~dimensions, is the case~$X_2=f(X_1)$ with $f$~strictly increasing. Then the only weighted points of~$\ens{m}$ are on the diagonal, in other words~$P(\mathbf{r},c)=m^{-1}\ind{r_1=r_2}$. On the other hand, if all components are independent, a symmetry argument gives $P(\mathbf{r},c)=m^{-d}$.
  \item In a practical setting $\binom{n}{m}$ will be far too large, so we will not be able to draw all sub-samples. We will use a random sub-sampling to obtain an approximation of~$\hat{P}_n$.
\end{itemize}

\subsection{Examples}
We propose a small example completely detailed. Table~\ref{exampledata} is a 4~sample in~$\mathbb{R}^2$ (each observation is identified by a lowercase letter), we choose~$n=4$ and~$m=3$. Table~\ref{examplecalcul} summarizes the computations.
\begin{table}
  \tbl{Example: data. $R_X$ (resp~$R_Y$) stands for the rank in~$X$ (resp.~$Y$)}
{\begin{tabular}{@{}ccccc}\toprule
Observation & $X$ & $Y$ & $R_X$ & $R_Y$\\
\colrule
 a &2.29&-0.97&4&1\\
 b &-1.2&-0.95&1&2\\
 c &-0.69&0.75&2&4\\
 d &-0.41&-0.12&3&3\\
   \botrule
  \end{tabular}}
\label{exampledata}
\end{table}

\begin{figure}
\begin{center}
\subfigure[Data]{
\resizebox*{5cm}{!}{\includegraphics{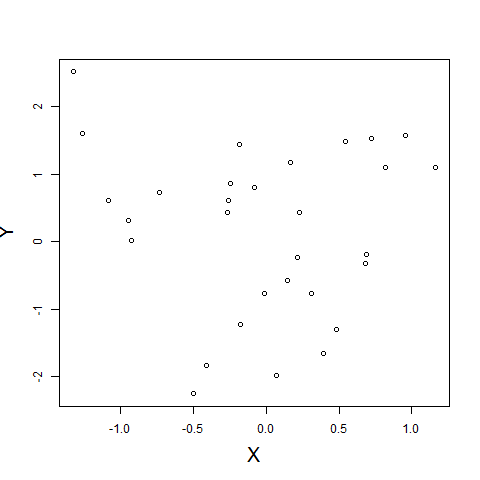}}}%
\subfigure[Estimation]{
\resizebox*{5cm}{!}{\includegraphics{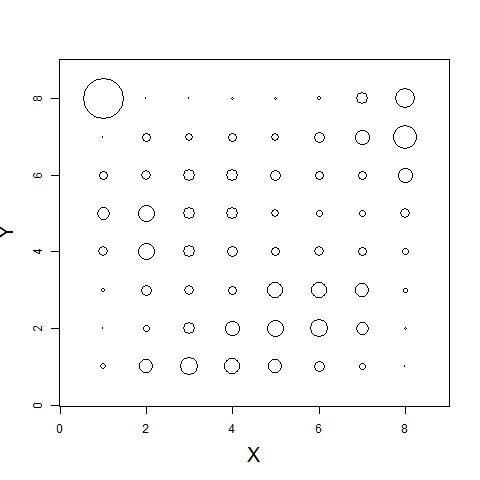}}}%
\caption{\label{plot30}Sample of 30~observations, such that~$Y=0.5 \cdot X^2+\epsilon$, $X$ and~$\epsilon$ are normally distributed, centred and reduced. The sub-sample size is~8. In the right-hand graph, the radius of each circle is proportional to~$\hat{P}_n(\mathbf{r})$.}
\end{center}
\end{figure}

\begin{table}
\tbl{Computation of~$\hat{P}_n(\mathbf{r})$ for the data of Table~\ref{exampledata}. The sub-samples are named:~$A=\{a,b,c\},B=\{a,b,d\},C=\{a,c,d\},D=\{b,c,d\}$. For example, $bD$ in the first cell of the first line means that observation~$b$ in sub-sample~$B$ is the first one in~$X$ and in~$Y$, and~1/12 is the value of~$\hat{P}_n((1,1))$, since~$12=3\times\binom{4}{3}$}
{\begin{tabular}{@{}ccrrr}\toprule
&&\multicolumn{3}{c}{Rank in~$X$}\\
&&1&2&3\\\colrule
\multirow{3}{*}{Rank in~$Y$}
& 1 & \{bD\}; 1/12 & $\varnothing$; 0/12 & \{aA,aB,aC\}; 3/12 \\
& 2 & \{bA,bB\}; 2/12 & \{dC\}; 1/12 & \{dD\}; 1/12 \\
& 3 & \{cC\}; 1/12 & \{cA,dB,cD\}; 3/12 & $\varnothing$; 0/12 \\
\botrule
\end{tabular}}\label{examplecalcul}
\end{table}
Figure~\ref{plot30} gives a more realistic example. We draw 30~observations, such that~$Y=0.5 \cdot X^2+\epsilon$, where~$X$ and~$\epsilon$ are normally distributed, centred and reduced. The sub-sample size is~8. In the right-hand graph, the radius of each circle is proportional to~$\hat{P}_n(\mathbf{r})$. One may remark that even in this small example, $m\times\binom{n}{m}\simeq \nbt{4}{7}$.

\section{Theoretical results}

\paragraph*{Why do we need an indirect study?} Addressing the behaviour of~$\hat{P}_n(\mathbf{r})$ is obviously a bit more difficult than addressing Kendall's~$\tau$, which is the case~$m=d=2$. Now, as far as we know, there is no simple formula for the probability distribution of Kendall's~$\tau$: neither ~\cite{wolfe1973nonparametric}, nor the~\cod{cor.test} function documentation (in the package~\cod{stats}) in \cod{R}, nor in the package~\cod{Kendall} of the same software, nor the \cod{SAS} documentation mention such a formula. There exist recursive ones in~\cite{BestGipps,valz1994exact}, which can be used only if the sample is small.
\paragraph*{Random sub-sampling} We noted previously that in many cases, $\binom{n}{m}$ is very large, so it is impossible to employ all the sub-samples to compute~$\hat{P}_n$: we need to select some of them randomly to obtain an approximation.
This kind of approximation has been addressed by~\cite{Blom76}: the additional variance is a function of a term which we do not know here ($\sigma_k$, using the notations of~\cite{Blom76}). On the other hand, we know that each cell of~$B_{\mathbf{r}}(X(S(\ens{m})))$ takes only 2~values: 0 and~$1/m$. This makes all convergence issues much easier.\\
Indeed, for a given sub-sample number~$b$, for each~$\mathbf{r}$, the number of sub-samples such that~$B(\mathbf{r},X(\ens{m}))=1/m$ is binomially distributed, with parameters~$b$ and~$\hat{P}_n(n,\mathbf{r})$. If~$b \hat{P}_n(n,\mathbf{r})>30$, this binomial distribution converges in distribution to a normal one, with standard deviation~$\left(b\hat{P}_n(n,\mathbf{r})\right)^{-1/2}$.\\
This normality allows of obtaining much information about any function of~$\hat{P}_n(n,\mathbf{r})$. In the worst case, the approximation errors are perfectly correlated, and the standard deviation of the global error is the sum of all the standard deviations.

\begin{theorem}\label{convponc} Assuming that~$\mathbf{X}=(X_1, \cdots, X_d)$ in~$\mathbb{R}^d$ has continuous marginals, if~$\mathbf{r}(m)/m \rightarrow \mathbf{x}$ when~$m$ tends to infinity, with~$c$ continuous in~$\mathbf{x} \in [0,1]^d$, then:
\begin{gather*}
  m^d\times P(\mathbf{r}(m),c) \rightarrow c(\mathbf{x})
.\end{gather*}
\end{theorem}
We can derive the following property:
\begin{proposition}\label{convT}  Assuming that~$\mathbf{X}=(X_1, \cdots, X_d)$ in~$\mathbb{R}^d$ has continuous marginals and continuous bounded copula:
\begin{gather*}
 m^{d} \times \sum_{\mathbf{r}}{\left(P(\mathbf{r},c) - P(\mathbf{r},c_0) \right)^2} \rightarrow
     \int_{[0,1]^d}{\left(c(\mathbf{x})-c_0(\mathbf{x})\right)^2 d\mathbf{x}}
.\end{gather*}
\end{proposition}
That is why, in the following, we will study
\begin{gather*}
 T = m^d \times \sum_{\mathbf{r}}{\left(\hat{P}_n(\mathbf{r}) - P(\mathbf{r},c_0) \right)^2}
.\end{gather*}
More precisely, we study~$\varloi{T}{c_0}$, where~$c_0$ stands for any distribution with globally independent components.
\begin{theorem}\label{varnul}
Assuming that~$\mathbf{X}=(X_1, \cdots, X_d)$ in~$\mathbb{R}^d$ has continuous marginals, if the components are globally independent and if the sample size~$n$ tends to infinity,
\begin{gather*}
\begin{array}{ll}
 \varloi{T}{c_0}\times n^2 \rightarrow &
     2S_2^d - 2m^4+\\
    &2(m-1)^4 \left(\frac{m^4-4m^3+6m^2-4m+S_2}{(m-1)^4}\right)^{d}+
    12(m-1)^2 \left(\frac{m^2-2m+S_2}{(m-1)^2}\right)^{d}+\\
    & 8(m-1) \left(\frac{m-S_2}{m-1}\right)^{d} +
    8(m-1)^3 \left(\frac{m^3-3m^2+3m+S_2}{(m-1)^3}\right)^{d}\\
 \esploi{T}{c_0} \times n \rightarrow &
     S_1^d -m^2 + (m-1)^2 \left(\frac{m^2-2m+S_1}{(m-1)^2}\right)^d+2(m-1)\left(\frac{m-S_1}{m-1}\right)^d
\end{array}
\end{gather*}
with
\begin{gather*}
 S_2=\frac{m^2\binom{4m-3}{2m-2}}{\left((2m-1)\binom{2m-2}{m-1}\right)^2}
 \ \ \ \ \ \ 
 S_1=\frac{m4^{m-1}}{(2m-1)\binom{2m-2}{m-1}}
.\end{gather*}
\end{theorem}
For large values of~$m$, using the first order Taylor expansion, we show
\begin{corollary}\label{varnulA}
\begin{gather*}
 \lim_{n\rightarrow \infty}{\left(\varloi{T}{c_0}\times n^2\right)} \simeq  2\left(\sqrt{\frac{\pi m}{8}}\right)^d +\left(7\sqrt{\frac{\pi}{2}}\right)d\sqrt{m}
\end{gather*}
\end{corollary}
For~$\esploi{T}{c_0}$, the same computation is a bit more difficult since we need a second order Taylor expansion, and much less useful since the mean value of a test statistic is not an issue.\\
Table~\ref{borderL} shows, with respect to~$m$, the value~$\max\left\{d| \left(\varloi{T}{c_0}\times n^2\right)/\left(\left(7\sqrt{\frac{\pi}{2}}\right)d\sqrt{m}\right)>1/2 \right\}$.
\begin{table}
  \tbl{For $m$, the value of~$d$ for which the linear part of variance becomes less than half the total variance}
{\begin{tabular}{@{}cc}\toprule
sub-sample size~$m$ & border value for~$d$\\
\colrule
10&  5\\
15&  4\\
20&  4\\
  \botrule
  \end{tabular}}
\label{borderL}
\end{table}

\section{Simulations}
We performed all these simulations using~\texttt{R}, linked with a~\texttt{DLL} written in~\texttt{C} (parallelized with~\texttt{Open-MP}~\citep{openmp}), for reasons of efficiency. The random number generator is described in~\citep{randomkit}.

\subsection{Convergence to the normal distribution}
\begin{figure}
\begin{center}
\subfigure[]{
\resizebox*{5cm}{!}{\includegraphics{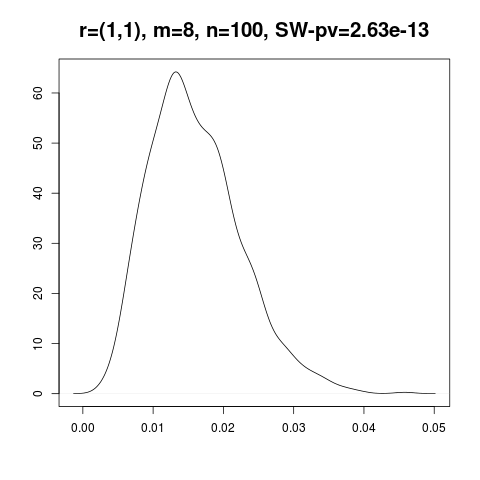}}}
\subfigure[]{
\resizebox*{5cm}{!}{\includegraphics{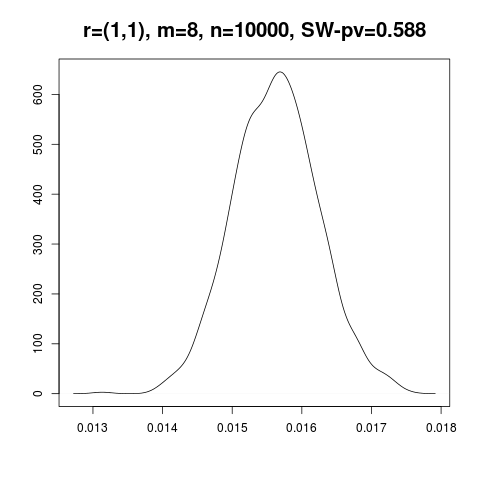}}}\\
\subfigure[]{
\resizebox*{5cm}{!}{\includegraphics{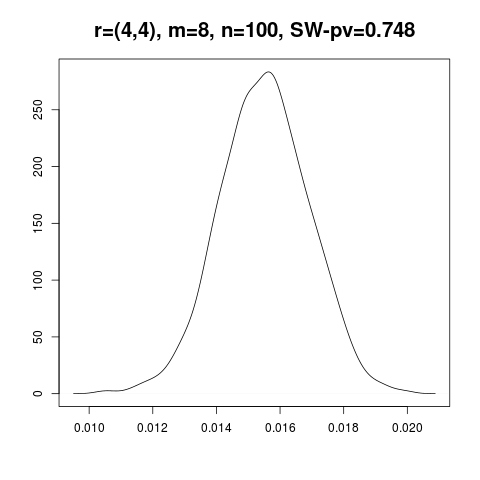}}}
\subfigure[]{
\resizebox*{5cm}{!}{\includegraphics{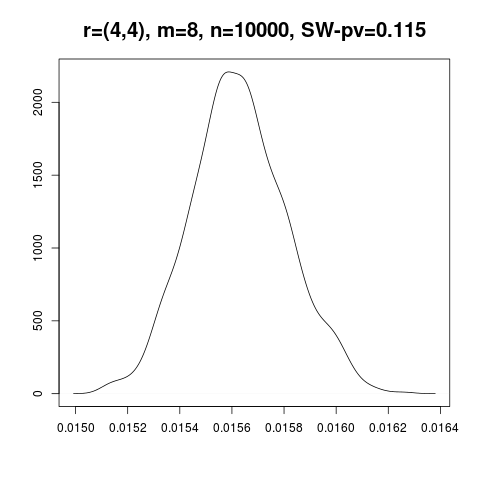}}}
\caption{\label{histos}Empirical densities of~$\hat{P}_n(\mathbf{r})$: $m$ and~$n$ have the usual meanings, SW-pv is the $p$-value of the Shapiro--Wilk test for 1000~samples.}
\end{center}
\end{figure}

We look into the convergence using the graphs in Figure~\ref{histos}. As the size of the sub-sample is~8, and the dimension~2, the vector of ranks~$\mathbf{r}=(1,1)$ is a corner, and so has a huge dispersion, and $(4,4)$ is the centre, with a small dispersion. The convergence is obviously better for places with small dispersion.

\subsection{Distance simulation}
We perform the following simulations, with independent components:
\begin{itemize}
  \item $ m \in \{10,15\} $
  \item $ d \in \{2,3,4,5\} $
  \item $ n \in \{100,10000\} $
\end{itemize}
and we cross all these possibilities, which makes 16~possibilities.\\
For all simulations, we choose a number of sub-samples equal to~$15n^d$, and a number of samples equal to~200.\\
Then we can assert that:
\begin{itemize}
  \item Convergence is good for means, probably because we only need the convergence of the variance of the $U$-statistic to the first term~$\sigma_1$.
  \item Convergence is slower for variances, probably because we additionally need the convergence to the normal distribution. The graphs in Figure~\ref{histos} show that this convergence is slow.
\end{itemize}
The complete results are given in Table~\ref{ressimsL2}.
\begin{table}
  \tbl{Means and variances: $m$ and~$v$ denote mean and variance, $\overline{x}$ denotes the ratio's true value / theoretical value for the parameter~$x$.}
{\begin{tabular}{@{}lllccrr}\toprule
$m$ &$d$ & $n$ & $m$ & $v$ & $\overline{m}$ & $\overline{v}$ \\\colrule
   10&    2&  100&4.8e-04&3.5e-08& 1.15& 0.05\\
   10&    3&  100&2.2e-04&1.7e-09& 1.17& 0.15\\
   10&    4&  100&7.5e-05&9.9e-11& 1.23& 0.59\\
   10&    5&  100&2.3e-05&9.8e-12& 1.28& 4.05\\
   15&    2&  100&3.8e-04&1.5e-08& 1.24& 0.09\\
   15&    3&  100&1.4e-04&3.7e-10& 1.30& 0.29\\
   15&    4&  100&3.7e-05&1.7e-11& 1.37& 1.83\\
   15&    5&  100&9.6e-06&1.7e-12& 1.49&25.42\\
   10&    2&10000&1.3e-05&7.8e-12& 3.16& 0.11\\
   10&    3&10000&3.0e-06&1.7e-13& 1.59& 0.15\\
   10&    4&10000&7.0e-07&4.4e-15& 1.15& 0.27\\
   10&    5&10000&1.9e-07&8.1e-17& 1.06& 0.34\\
   15&    2&10000&5.9e-06&1.4e-12& 1.90& 0.08\\
   15&    3&10000&1.2e-06&2.0e-14& 1.18& 0.16\\
   15&    4&10000&2.8e-07&3.1e-16& 1.05& 0.34\\
   15&    5&10000&6.6e-08&4.6e-18& 1.02& 0.68\\
   \botrule
  \end{tabular}}
\label{ressimsL2}
\end{table}

\subsection{Use for independence test}\label{testtest}
We propose here to use independence testing as a measure of the accuracy of our estimation technique.\\
In order to study the behaviour of this testing method with high-dimensional data, we will use as a test case a fixed dependence between~$X_1$ and~$X_2$, with all other~$X_i$ being independent. The fixed dependence is~$X_2=0.5\cdot X_1^p + \epsilon$, $p=1,2$ (if~$p=1$, the dependence is monotonic, and non-monotonic otherwise). The random variable~$\epsilon$, plus all~$X_i$ except~$X_2$, have a standard Gaussian distribution. The level of the test is~5\%.\\
The theory developed previously does not allow explaining the empirical observations. Indeed:
\begin{itemize}
  \item The behaviour of the distance with the uniform case is not within the scope of the theoretical study;
  \item Furthermore, the dissimilarity used for the test is not the $L^2$~distance but the Kullback--Leibler divergence.
\end{itemize}
We preferred the Kullback--Leibler divergence because it gives better results when, in each set of probabilities, there are very heterogeneous values. The Kullback--Leibler divergence is computed as follows:
\begin{gather*}
  \hat{T}(c_0)=
  \sum_{\mathbf{r}\in \ens{m}^d}
  {\hat{P}_n(\mathbf{r}) \log\left(\frac{\hat{P}_n(\mathbf{r})}{P(\mathbf{r},c_0)}\right)}
.\end{gather*}
To build a test, one just has to simulate many samples with independent components, which gives the distribution of the test statistic.
\subsubsection{Independence test methods} As we focus here on testing independence, we improve the state of the art proposed in the Introduction.\\
Some methods more specific to independence testing have been recently developed. The method in~\cite{Panchenko2005176} uses $V$-statistics to estimate a distance between probability distributions. \citep{Bakirov20061742} uses a measure of association determined by inter-point distances. \cite{KalLed1999} uses Fourier decomposition to summarize the density of the rank of the observations. \cite{2009arXiv0908.2794G} focuses on the case of two variables: it is based on the size of the longest increasing subsequence of the permutation which maps the ranks of the~$X$ observations to the ranks of the~$Y$ observations.\\
A natural tool for comparing these heterogeneous methods is a simulation study. \cite{berggof} is a comprehensive simulation study. It focuses on Deheuvels or Rosenblatt type tests. The main conclusion of that paper is that the original Deheuvels test is among the best of this family. However, \cite{2007arXiv0709.3860C} states that Deheuvels tests are poor, in terms of efficiency, when the components of the multivariate random variable are not positively dependent (A random variable~$(X, Y)$ is positively dependent if and only if~$\forall (x,y) \in \mathbb{R}^2 \ \proba{X\leq x,Y\leq y} \geq \proba{X\leq x}\proba{Y\leq y}$). In this case the half-spaces we use to detect the deviation from uniformity are not simultaneously overcharged (or simultaneously undercharged), so the test statistic remains small. On the other hand, as far as we know, there is no reason to think these tests are affected by high dimensionality.\\
It is the opposite for kernel goodness-of-fit tests. The simulation study included in~\cite{sca} shows that they are able to detect deviations with a very complicated form. About the impact of dimensionality, we conjecture that these tests suffer from the well known \gui{curse of dimensionality}.\\
Based on \cite{berggof}, we propose a partial simulation study, including: the new test using sub-sampling, the original Deheuvels test~\citep{deheuvels1}, and a test using kernel estimation, studied in~\cite{sca}, which is one of the most recent studies in this vein.
\subsubsection{Implementation} We used~1000~samples for the simulation of the test statistic distribution under the independence hypothesis, and 300~for the test power evaluation for the alternate hypothesis.
\begin{itemize}
  \item The Deheuvels test implementation was that of the \texttt{R}~package~\texttt{copula}.
  \item The kernel goodness-of-fit test was implemented using the \texttt{R}~package~\texttt{ks}. We succeeded in reproducing the results summarized in~\citep{sca}, using a slightly different bandwidth choice method. Scaillet used \gui{Scott's rule of thumb}, modified by an empirically chosen factor, to maximize the test's power. This factor was often equal to~0.5. In the package~\texttt{ks}, the fastest bandwidth choice method implemented was~\texttt{Hpi}, and computation time was a constraint for this study, so we used it. This led us to choose~0.33 in place of~0.5 to maximize the test's power. We could not study dimensions greater than~3, because of computer limitations.
  \item For the sub-sampling method, the number of sub-samples was~$m=10^5$, and the sub-sample size~8.
\end{itemize}
\subsubsection{Results} Let us look at the main features of the results, summarized in Table~\ref{puissances_dim}.
\begin{itemize}
  \item For the Deheuvels test, we confirmed the poor power for non-monotonic dependences, and a low impact of dimensionality.
  \item We almost confirmed our conjectures for the kernel goodness-of-fit test: the form of the dependence does not change markedly the power of the test, and the impact of the dimension~3 is important. We can not say more about the impact of dimensionality.
  \item About the sub-sampling test, we can assert that it suffices to increase linearly the number of observations to maintain the same test power (with a fixed 5\%~level).
\end{itemize}
\begin{table}
  \tbl{Power w.r.t. dimension: ``Deheuvels'' stands for the Deheuvels test of independence, and ``Kernel'' for a GoF test based on kernel density estimation. Some values are missing because of memory limitations.}
{\begin{tabular}{@{}cccccc}\toprule
Dependence & Dimension & Sample size& sub-sampling & Deheuvels & Kernel \\
\colrule
$X_2=0.5X_1^2+\epsilon$ & 2 & 30 & 0.44 & 0.13 & 0.20 \\
$X_2=0.5X_1^2+\epsilon$ & 3 & 45 & 0.41 & 0.06 & 0.15 \\
$X_2=0.5X_1^2+\epsilon$ & 4 & 60 & 0.39 & 0.09 & \\
$X_2=0.5X_1^2+\epsilon$ & 6 & 90 & 0.33 & 0.08 & \\\colrule
$X_2=0.5X_1+\epsilon  $ & 2 & 30 & 0.54 & 0.63 & 0.35 \\
$X_2=0.5X_1+\epsilon  $ & 3 & 45 & 0.54 & 0.36 & 0.20 \\
$X_2=0.5X_1+\epsilon  $ & 4 & 60 & 0.53 & 0.28 & \\
$X_2=0.5X_1+\epsilon  $ & 6 & 90 & 0.48 & 0.22 & \\
   \botrule
  \end{tabular}}
\label{puissances_dim}
\end{table}
The results show a specific behaviour of the sub-sampling method, regarding average dimensions.
\subsubsection{Other simulations} It does not seem that the linear increase of the number of observations is due to the number of dependent variables~(2). We simulate the following model:
\begin{gather*}
\begin{array}{lll}
  V \sim \mathcal{LN}(0,a)&&\\
  X_i \sim \mathcal{N}(0,V) & \textrm{ if } & i \leq d_d\\
  X_i \sim \mathcal{N}(0,1) & \textrm{ if } & d_d< i \leq d
\end{array}
\end{gather*}
where~$\mathcal{LN}$ denotes the lognormal distribution. In other words, the first~$d_d$ components have the same random volatility. If $a$~is large, the dependence is strong, and null if~$a=0$. So, whether~$d_d$ is equal to~2 or~3, we obtain the same result as before.

\subsection{Use for regression}
\begin{figure}
\begin{center}
\subfigure[\textbf{True} model]{
\resizebox*{5cm}{!}{\includegraphics{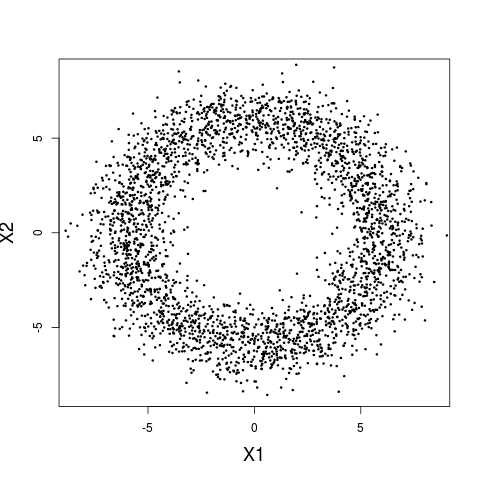}}}
\subfigure[\textbf{Estimated} model]{
\resizebox*{5cm}{!}{\includegraphics{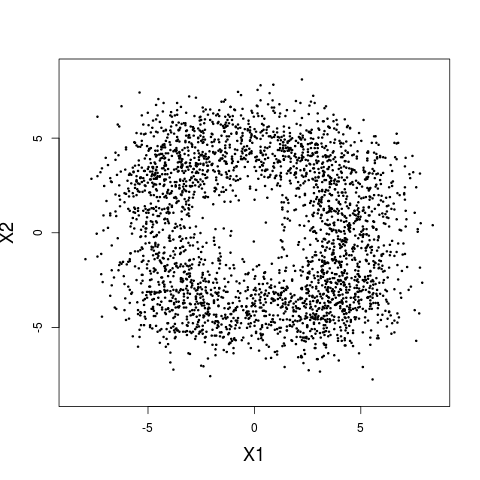}}}
\caption{\label{boules}Example of estimation. True model: $\mathbf{X}=a\mathbf{N}_1/|\mathbf{N}_1|+\mathbf{N}_2$, with~$\mathbf{N}_i \sim \mathcal{N}(0,I_5)$, and~$a=6$. Estimation on 3000~observations. Conditioning: $X_3=X_4=X_5=0$, which gives the same ``sphere'' in dimension~2.}
\end{center}
\end{figure}

Using Theorem~\ref{convponc}, it is possible to use the estimation~$\hat{P}_n$ to build a regression model.\\
The first step is to transform a law on~$\ens{m}^d$ into a law on~$[0,1]^d$. We could use~$P(\ent{m\mathbf{x}}+1,c)$, as in the proof of property~\ref{convT}, this is a convergent interpolation, but not very smooth. We preferred to use~$\sum_{\mathbf{r}}{\hat{P}_n(\mathbf{r})\prod_{l=1}^d{b_{m,r_l}(x_l)}}$: it is convergent too (it is a Bernstein approximation), much smoother, and easy to implement (it comes down to simulating $\beta$~distributions).\\
The second step is to use marginal densities to transform the copula into the original density: we used the usual kernel density estimation.\\
To show the efficiency of this estimation method, we use a simple model. The results are shown in Figure~\ref{boules}. This model is convenient because of its simplicity and scalablity (it is easy to change the dimension of the model). Furthermore, it is difficult to estimate, because any line not too far from~0 is first in a sparse region, then in a dense one, and then sparse, dense, and sparse again.

\section{Ackowledgements}
We are pleased to thank Prs~Bertail and Clemen\c{c}on for their help, comments, and remarks.

\appendices
\section{Proofs of Theorem~\ref{convponc} and property~\ref{convT}}
\subsection{Theorem~\ref{convponc}}
In this part, we assume~$\mathbf{X}$ has uniform marginal distributions on~$[0,1]$, in other words its distribution is a copula, denoted~$c$ in the following.\\
We compute the probability that the rank~$\mathbf{r}(m)$ is reached by the last observation. Using symmetry, the probability that any of the observations reaches the rank~$\mathbf{r}(m)$~is $m$~times bigger, but one has afterwards to divide by~$m$ to get~$P$.\\
In order for the last observation to reach rank~$\mathbf{r}(m)$, one needs its first coordinate to be between~$X_{1,[r_1(m)-1,m-1]}$ and~$X_{1,[r_1(m),m-1]}$, where~$X_{1,[r_1(m)-1,m-1]}$ denotes the value with rank~$r_1(m)-1$ amongst~$m-1$, for the first coordinate. It is the same for the other coordinates, so
\begin{gather}\label{trou}
  P(\mathbf{r}(m),c) = \esp{\int_{ \bigcap{[X_{l,[r_l(m)-1,m-1]};X_{l,[r_l(m),m-1]}]} }{c(\mathbf{u})d\mathbf{u}}}
.\end{gather}
Then
\begin{gather*}
    \esp{ \prod{\left(X_{l,[r_l(m),m-1]}-X_{l,[r_l(m)-1,m-1]}\right)}} \times 
    \min_{\bigcap_{l=1}^d{ \left[ X_{l,[r_l(m)-1,m-1]},X_{l,[r_l(m),m-1}]\right] }}{c(\mathbf{x})}
    \\ \leq P(\mathbf{r}(m),c) \leq \\
    \esp{ \prod{\left(X_{l,[r_l(m),m-1]}-X_{l,[r_l(m)-1,m-1]}\right)}} \times 
    \max_{\bigcap_{l=1}^d{ \left[ X_{l,[r_l(m)-1,m-1]},X_{l,[r_l(m),m-1}]\right] }}{c(\mathbf{x})}
.\end{gather*}
On the other hand, both $\left(X_{1,[r_1(m)-1,m-1]}, \cdots, X_{d,[r_d(m)-1,m-1]}\right)$ and~$\left(X_{1,[r_1(m),m-1]}, \cdots, X_{d,[r_d(m),m-1]}\right)$ converge in probability to~$r(m)/m$, and~$c$ is continuous, so
\begin{gather*}
  P(\mathbf{r}(m),c) \rightarrow 
    \esp{ \prod{\left(X_{l,[r_l(m),m-1]}-X_{l,[r_l(m)-1,m-1]}\right)}} \times c(r(m)/m)
.\end{gather*}
We now have to study the expectation. One notes~$\mathbf{x}$ and~$\mathbf{\epsilon}$ two vectors in~$|0,1]^d$, and~$\mathbf{r}(m)$ a vector of ranks, so in~$\ens{m}^d$. We want to compute
\begin{gather*}
  \proba{ \forall l \ X_{l,[r(m)_l-1,m-1]}<x_l \cap X_{l,[r_l(m),m-1]}>x_l+\epsilon_l}=\\
    \proba{ \forall l \ X_{l,[r_l-1,m-1]}<x_l } \times \\
    \proba{\textrm{no observation in }\bigcup_{l=1}^d{[x_l,x_l+\epsilon_l]}}
,\end{gather*}
because the two regions are to have empty intersection. We have\\
\begin{gather*}
  \proba{\textrm{no observation in }\bigcup_{l=1}^d{[x_l,x_l+\epsilon_l]}}=
  \left(1-C\left(\bigcup_{l=1}^d{[x_l,x_l+\epsilon_l]}\right)\right)^{m-1}
,\end{gather*}
and now
\begin{gather*}
  C\left(\bigcup_{l=1}^d{[x_l,x_l+\epsilon_l]}\right) \leq \sum_{l=1}^d{\epsilon_l}
.\end{gather*}
This is true for any value of~$\mathbf{x}$, so
\begin{gather*}
  \proba{\forall l,\ X_{l,[r_l(m),m-1]}-X_{l,[r_l(m)-1,m-1]}>\epsilon_l} \geq \left(1-\sum_{l=1}^d{\epsilon_l}\right)^{m-1}
,\end{gather*}
and, using integration by parts,
\begin{gather*}
  \esp{ \prod{\left(X_{l,[r_l(m),m-1]}-X_{l,[r_l(m)-1,m-1]}\right)} }=\\
  \int_{\mathbf{R}^{+d}}{\proba{\forall l,\ X_{l,[r_l(m),m-1]}-X_{l,[r_l(m)-1,m-1]}>\epsilon_l} d\epsilon_1 \cdots \epsilon_d} \geq\\
  \int_{\mathbf{R}^{+d}}{\left(1-\sum_{l=1}^d{\epsilon_l}\right)^{m-1}\epsilon_1 \cdots \epsilon_d}
.\end{gather*}
We note~$s=\sum_{l=1}^d{\epsilon_l}$ and
\begin{gather*}
  \esp{ \prod{\left(X_{l,[r_l(m),m-1]}-X_{l,[r_l(m)-1,m-1]}\right)} } \geq \\
  \int_{\mathbf{R}^{+}}{\frac{s^{d-1}}{(d-1)!}\left(1-s\right)^{m-1} ds}=
  \frac{(m-1)!}{(m+d-1)!} \rightarrow (m-1)^{-d}
,\end{gather*}
so then
\begin{gather*}
  \liminf_{m\rightarrow \infty} \left(m^d\times P(\mathbf{r}(m),c) )\right) \geq c(r(m)/m) \rightarrow c(\mathbf{x})
,\end{gather*}
and now
\begin{gather*}
  \sum_{\mathbf{r}}{P(\mathbf{r},c)}=1\ \ \ 
  \int_{[0,1]^d}{c(\mathbf{x})d\mathbf{x}}=1
,\end{gather*}
which gives the conclusion.

\subsection{Property~\ref{convT}}
We obviously have
\begin{gather*}
  m^d \times \sum_{\mathbf{r}}{\left(P(\mathbf{r},c) - P(\mathbf{r},c_0) \right)^2} =\\
    m^{2d} \times\int_{[0,1]^d}{\left(P(\ent{m\mathbf{x}}+1,c)-P(\ent{m\mathbf{x}}+1,c_0)\right)^2 d\mathbf{x}}=\\
  \int_{[0,1]^d}{\left(m^d P(\ent{m\mathbf{x}}+1,c)-m^d P(\ent{m\mathbf{x}}+1,c_0)\right)^2 d\mathbf{x}}
.\end{gather*}
Now~$(\ent{m\mathbf{x}}+1)/m \rightarrow \mathbf{x}$ when~$m$ tends to infinity, so~$m^d P(\ent{m\mathbf{x}}+1,c)\rightarrow c(\mathbf{x})$. We now have to prove that this simple convergence is dominated to prove the integral is convergent.\\
We use again equality~\ref{trou}. We have
\begin{gather*}
  C\left(\bigcup_{l=1}^d{[x_l,x_l+\epsilon_l]}\right) \geq 
  \max_{l=1}^d{\epsilon_l}
,\end{gather*}
which gives
\begin{gather*}
  \proba{\forall l,\ X_{l,[r_l(m),m-1]}-X_{l,[r_l(m)-1,m-1]}>\epsilon_l} \leq \left(1-\max_{l=1}^d{\epsilon_l}\right)^{m-1} \leq e^{-(m-1)\max_{l=1}^d{\epsilon_l}}
.\end{gather*}
This tail distribution is the one of the $d$-uple of random variables~$Z_1, \cdots, Z_d$, with~$Z_1= \cdots=Z_d=Z$, and~$Z$ exponentially distributed with parameter~$m$. Then
\begin{gather*}
 \esp{ \prod{\left(X_{l,[r_l(m),m-1]}-X_{l,[r_l(m)-1,m-1]}\right)} } \leq \esp{Z^d}=d!(m-1)^{-d}
,\end{gather*}
and
\begin{gather*}
 P(\mathbf{r},c) \leq d!(m-1)^{-d} \max_{[0,1]^d}(c)
,\end{gather*}
and so
\begin{gather*}
\forall m,\ \int_{[0,1]^d}{\left(m^d P(\ent{m\mathbf{x}}+1,c)-m^d P(\ent{m\mathbf{x}}+1,c_0)\right)^2 d\mathbf{x}} \leq d!\max_{[0,1]^d}(c)^2
,\end{gather*}
showing that these integrals are convergent.

\section{Proof of Theorem~\ref{varnul}}
\subsection{The \textit{U}-statistics: reminders and notations}
Let $h$~be a measurable function, symmetric in its $n$~arguments. Then if we have a sample~$X_1, \cdots, X_N$ with~$n>m$, we define the $U$-statistic~$U_m$:
\begin{gather*}
  U_n = \frac{1}{\binom{n}{m}} \sum_{S \subset \ens{N}}^{\card{S}=m}
  {h\left(X_{S(1)},\cdots,X_{S(m)}\right)}
.\end{gather*}
where~$S(i)$ denotes the \ordinal{i} element of~$S$.\\
Furthermore, we define
\begin{gather*}
  h_c = \esp{h\left(x_1, \cdots, x_c,X_{c+1}, \cdots, X_m\right)}
,\end{gather*}
and
\begin{gather*}
  \sigma_c = \var{h_c\left(X_1, \cdots, X_c\right)}
.\end{gather*}
Then~$\esp{U_m}=h_0$ and, when~$n \rightarrow +\infty$, $n(U_m-h_0)$ converges in distribution to~$\mathcal{N}\left(0,m^2\sigma_1\right)$ if~$\sigma_1 \neq 0$.\\
Furthermore, if we have another $U$-statistic~$V_n$ defined by a kernel~$g$, we may also define~$g_c$ and~$\sigma_{c,c}$:
\begin{gather*}
  \sigma_{c,c} = \cov{h_c\left(X_1, \cdots, X_c\right),g_c\left(X_1, \cdots, X_c\right)}
.\end{gather*}
The covariance between~$U_n$ and~$V_n$ converges to~$\sigma_{1,1}$ when~$N \rightarrow +\infty$.\\
In the following, for each~$\mathbf{r}$, we have a $U$-statistic~$\hat{P}_n(\mathbf{r})$, whose normal convergence we will use.\\
As we are in a slightly special case of $U$-statistics (we are only interested in the case~$c=1$, but we study a large number of $U$-statistics at the same time), one has to adapt the notations. We note:
\begin{gather*}
  h(\mathbf{r},x_1)=\esp{B_\mathbf{r}\left(x_1,X_2\cdots,X_n\right)}\\
  \sigma(\mathbf{r},\mathbf{s})=\cov{h(\mathbf{r},X_1),h(\mathbf{s},X_1)}
.\end{gather*}
So we obtain
\begin{proposition}\label{hajekloc}
If~$n \rightarrow \infty$, one has
\begin{itemize}
  \item $n\left(\hat{P}_n(\mathbf{r}-P(\mathbf{r},c))\right)$ converges in distribution to~$\mathcal{N}\left(0,m^2\sigma(\mathbf{r},\mathbf{r})\right)$; if~$\sigma(\mathbf{r},\mathbf{r}) \neq 0$,
  \item $n\cov{\hat{P}_n(\mathbf{r}),\hat{P}_n(\mathbf{s})} \rightarrow m^2\sigma(\mathbf{r},\mathbf{s})$.
\end{itemize}
\end{proposition}
This Central Limit Theorem allows using the following computation:
\begin{proposition}\label{varD}
Let $\mathbf{X}$ be a vector such that$\mathbf{X} \leadsto \mathcal{N}(0,V)$, where the coefficients of~$V$ are denoted~$\sigma(\mathbf{r},\mathbf{s})$. We study~$D=\sum_{\mathbf{r}}{X_{\mathbf{r}}^2}$. Then
\begin{gather*}
  \esp{D}=\sum_{\mathbf{r}}{\sigma(\mathbf{r},\mathbf{r})}\ \ \ \ \ 
  \var{D}=2\sum_{\mathbf{r},\mathbf{s}}{ \sigma(\mathbf{r},\mathbf{s})^2 }
.\end{gather*}
\end{proposition}

\subsection{Calculation of the covariances}
\begin{figure}
\begin{center}
\resizebox*{10cm}{!}{\includegraphics{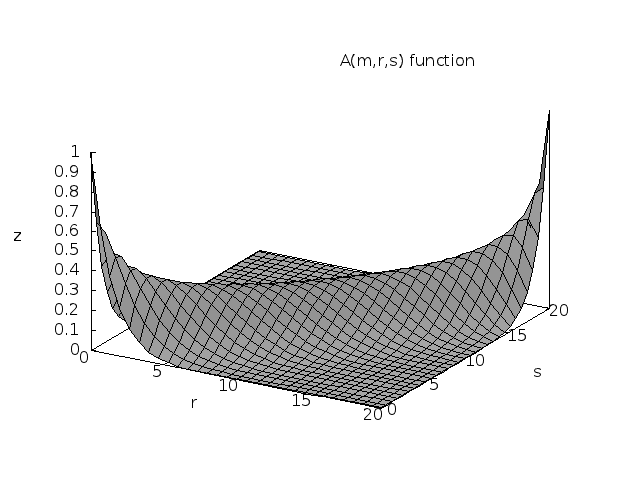}}
\caption{\label{amrs} Plot of $\mathcal{A}(m,r,s)$}%
\end{center}
\end{figure}

We calculate in the same way the variances and covariances. We consider the first observation, its vector of ranks, and a given rank~$\mathbf{r}$. There are 3~cases:
\begin{itemize}
  \item Other observations are such that the vector of ranks of the first observation is~$\mathbf{r}$,
  \item Other observations are such that the vector of ranks of the first observation is equal to~$\mathbf{r}$ for some dimensions,
  \item Other observations are such that the vector of ranks of the first observation is different of~$\mathbf{r}$ for all dimensions.
\end{itemize}
We know the probability that the \ordinal{l} coordinate of the first observation reaches rank~$r_l$:
\begin{gather*}
  P=\binom{m-1}{r_l-1} x_{1,l}^{r_l-1}(1-x_{1,l})^{m-r_l}=b_{m-1,r_l-1}(x_{1,l})
,\end{gather*}
where $b_{m-1,r_l-1}$ is a Bernstein polynomial, with well known properties, for example
\begin{gather*}
  \int_0^1{b_{m,r}(x)dx}=\frac{1}{m+1}
.\end{gather*}
We use this to calculate the probability of each one of the three cases:
\begin{gather*}
  \proba{B_\mathbf{r}\left(x_1,X_2\cdots,X_n\right)=\frac{1}{m}}=\\
  1 \times \prod_{l=1}^d{b_{m-1,r_l-1}(x_{1,l})}+
  0\times\left(1-\prod_{l=1}^d{\left(1-b_{m-1,r_l-1}(x_{1,l})\right)}-\prod_{l=1}^d{b_{m-1,r_l-1}(x_{1,l})}\right)+\\
  \frac{m-1}{(m-1)^d}\times\prod_{l=1}^d{\left(1-b_{m-1,r_l-1}(x_{1,l})\right)}
.\end{gather*}
One can remark that integrating this probability over~$\mathbf{x}_1$ gives back the unconditional probability
\begin{gather*}
  \proba{B_\mathbf{r}\left(X_1,X_2\cdots,X_n\right)=\frac{1}{m}}=\\
  1 \times \int_{[0,1]^d}{\prod_{l=1}^d{b_{m-1,r_l-1}(x_{1,l})}d\mathbf{x}_1}+
  \frac{1}{(m-1)^{d-1}}\times\int_{[0,1]^d}{\prod_{l=1}^d{\left(1-b_{m-1,r_l-1}(x_{1,l})\right)}d\mathbf{x}_1}=\\
  1\times\frac{1}{m^d}+\frac{1}{(m-1)^{d-1}}\times\left(1-\frac{1}{m}\right)^d=\frac{1}{m^{d-1}}
.\end{gather*}
The conditional mean is
\begin{gather*}
  h(\mathbf{r},\mathbf{x}_1)=\frac{1}{m}\esp{B_\mathbf{r}\left(\mathbf{x}_1,X_2\cdots,X_n\right)}=\\
  \frac{1}{m} \times
  \left(
    \prod_{l=1}^d{b_{m-1,r_l-1}(x_{1,l})}+
    \frac{1}{m(m-1)^{d-1}}\times\prod_{l=1}^d{\left(1-b_{m-1,r_l-1}(x_{1,l})\right)}
  \right)
.\end{gather*}
So, we have to calculate
\begin{gather*}
  \sigma(\mathbf{r},\mathbf{s})=\esp{
  \left[h(\mathbf{r},X_1) - \esp{h(\mathbf{r},X_1)} \right] \times  \left[h(\mathbf{s},X_1) - \esp{h(\mathbf{s},X_1)} \right]
    } = \\
  \int_{[0,1]^d}{
  \begin{array}{c}
  \left[
  \frac{1}{m}  \prod_{l=1}^d{b_{m-1,r_l-1}(x_{1,l})}+
  \frac{1}{m(m-1)^{d-1}}\prod_{l=1}^d{\left(1-b_{m-1,r_l-1}(x_{1,l})\right)}
  -\frac{1}{m^d}
  \right] \times \\
  \left[
  \frac{1}{m}  \prod_{l=1}^d{b_{m-1,s_l-1}(x_{1,l})}+
  \frac{1}{m(m-1)^{d-1}}\prod_{l=1}^d{\left(1-b_{m-1,s_l-1}(x_{1,l})\right)}
  -\frac{1}{m^d}
  \right] 
  \end{array} d\mathbf{x}_1}
.\end{gather*}
For these covariance computations, we will calculate expressions such as
\begin{gather*}
  \int_0^1{b_{m,r}(x)b_{m,s}(x)dx}=
  \frac{\binom{m}{r}\binom{m}{s}}{\binom{2m}{r+s}}\times\frac{1}{2m+1}
.\end{gather*}
We note
\begin{gather*}
  \mathcal{A}(m,r,s)=\frac{\binom{m}{r}\binom{m}{s}}{\binom{2m}{r+s}}
.\end{gather*}
One may look at the graph of the function~$\mathcal{A}(m,r,s)$, drawn in Figure~\ref{amrs}.
We will need to know some sums involving~$\mathcal{A}(m,r,s)$, they are proved in~\ref{calcomb}.
\begin{gather*}
 \sum_{r}{\mathcal{A}(m,r,s)}=\frac{2m+1}{m+1},\ \ \ 
 \sum_{r}{\mathcal{A}(m,r,r)}=\frac{4^m}{\binom{2m}{m}},\ \ \ 
 \sum_{r,s}{\left(\mathcal{A}(m,r,s)\right)^2}=\frac{\binom{4m+1}{2m}}{(\binom{2m}{m})^2}
\end{gather*}
We multiply two sums of three terms. We denote by~$D_{i,j}$ the product of terms numbered~$i$ in the first sum and~$j$ in the second one. 
\begin{gather*}
\begin{array}{lll}
  D_{1,1}&=&\int_{[0,1]^d}{\left[
    \frac{1}{m} \times \prod_{l=1}^d{b_{m-1,r_l-1}(x_{1,l})b_{m-1,s_l-1}(x_{1,l})}
  \right]dx}\\
  &=&
  \frac{1}{m^2(2m-1)^{d}} \times \prod_{l=1}^d{\mathcal{A}(m-1,r_l-1,s_l-1)}\\
  D_{2,2}&=&\int_{[0,1]^d}{\left[
    \frac{1}{m^2(m-1)^{2d-2}} \times\prod_{l=1}^d{\left((1-b_{m-1,r_l-1}(x_{1,l}))(1-b_{m-1,s_l-1}(x_{1,l}))\right)}
  \right] dx}\\
  &=&
  \frac{1}{m^2(m-1)^{2d-2}}\times\prod_{l=1}^d{\left(1-\frac{2}{m}+\frac{\mathcal{A}(m-1,r_l-1,s_l-1)}{2m-1}\right)}\\
  D_{3,3}&=&\int_{[0,1]^d}{\left[\frac{1}{m^{d}}\right]^2 dx}\\
  &=&\frac{1}{m^{2d}}\\
  D_{1,2}&=&\int_{[0,1]^d}{
    \frac{1}{m^2(m-1)^{d-1}} \times
    \prod_{l=1}^d{\left(b_{m-1,r_l-1}(x_{1,l})-b_{m-1,r_l-1}(x_{1,l})b_{m-1,s_l-1}(x_{1,l})\right)}
  dx}\\
  &=&\frac{1}{m^2(m-1)^{d-1}} \times \prod_{l=1}^d{\left(\frac{1}{m}-\frac{\mathcal{A}(m-1,r_l-1,s_l-1)}{2m-1}\right)}\\
  D_{1,3}&=&\int_{[0,1]^d}{\left[
    \frac{1}{m} \times \prod_{l=1}^d{b_{m-1,r_l-1}(x_{1,l})} \times \frac{1}{m^d}
  \right]dx}\\
   &=&\frac{1}{m^{2d+1}}\\
  D_{2,3}&=&2\int_{[0,1]^d}{\left[
    \frac{1}{m^{d+1}(m-1)^{d-1}} \times\prod_{l=1}^d{\left(1-b_{m-1,r_l-1}(x_{1,l})\right)}
  \right]dx}\\
  &=&\frac{m-1}{m^{2d+1}}
\end{array}
\end{gather*}
It is clear that~$D_{1,2}=D_{2,1}$, $D_{1,3}=D_{3,1}$ and~$D_{2,3}=D_{3,2}$, so we write simply~$2D_{1,2}$, etc.\\
We remark that
\begin{gather*}
  D_{3,3}-2D_{1,3}-2D_{2,3}=-m^{-2d}.
\end{gather*}
The sums built from the three remaining terms are sums of products, we transform them easily into products of sums. For example
\begin{gather*}
  \sum_{\mathbf{r},\mathbf{s}}{D_{1,1}^2} = 
  \frac{1}{m^4(2m-1)^{2d}} \times
    \sum_{\mathbf{r},\mathbf{s}}\left(\prod_{l=1}^d{\mathcal{A}^2(m-1,r_l-1,s_l-1)}\right)=\\
  \frac{1}{m^4(2m-1)^{2d}} \times
    \left(\sum_{r,s}\mathcal{A}^2(m-1,r-1,s-1)\right)^d
.\end{gather*}
All of these sums (3~sums of degree~1 with~$\mathbf{r}=\mathbf{s}$, 3~sums of degree~1, 3 sums of squares, and 3~sums of double products) are calculated in~\ref{somD}. Using these sums, the proof of Theorem~\ref{varnul} is obvious.

\section{Tools for the proof of Theorem~\ref{varnul}}
\subsection{Combinatorial computations}\label{calcomb}
We need to compute some sums involving~$\mathcal{A}(m,r,s)$. We first remark:
\begin{gather*}
  \mathcal{A}(m,r,s)=\frac{\binom{m}{r}\binom{m}{s}}{\binom{2m}{r+s}}=\frac{\binom{r+s}{r}\binom{2m-r-s}{m-r}}{\binom{2m}{m}}
.\end{gather*}
We show
\begin{proposition}{If~$m>1$ and~$0\leq s \leq m$:
\begin{gather*}
 \sum_{r}{\binom{r+s}{r}\binom{2m-r-s}{m-r}}=\binom{2m+1}{m}\\
 \sum_{r}{\binom{2r}{r}\binom{2m-2r}{m-r}}=4^m\\
 \sum_{r,s}{\left(\binom{r+s}{r}\binom{2m-r-s}{m-r}\right)^2}=\binom{4m+1}{2m}
.\end{gather*}
}
\end{proposition}
These sums are very similar to convolutions, so it is interesting to use generating functions~\citep{gfology}. They are quite simple for the first series:
\begin{gather*}
  \sum_n{\binom{n+k}{n}x^n}=\frac{1}{(1-x)^{k+1}}\\
  \sum_n{\binom{2n}{n}x^n}=\frac{1}{\sqrt{1-4x}}
.\end{gather*}
It is also possible to demonstrate the first identity using combinatorial arguments, counting the number of paths joining two opposite corners of a rectangle with sides~$m$ and~$m+1$. The last identity is more difficult: one needs to use the powerful tools developed in~\cite{aeqb}. As the sum is over 2~variables, one needs to use the package~\cod{multisum}~\citep{multisum}. The code proving identity is
\begin{verbatim}
Get["C:\Users\Jerome\Desktop\Celine\MultiSum.m"]
FindRecurrence[ (Binomial[r+s,r]*Binomial[2*m-r-s,m-r])^2, m, {r,s}, 4 ]
SumCertificate[%]
CheckRecurrence[ %, Binomial[4*m+1,2*m] ]
\end{verbatim}
For asymptotic expressions, we will use
\begin{gather*}
 \binom{2m}{m}=\frac{4^m}{\sqrt{\pi m}}\left(1-\frac{c_m}{m}\right)
,\end{gather*}
where~$\forall m, 1/9<c_m<1/8$.

\subsection{Sums of covariance terms}\label{somD}
We summarize here the sums of terms~$D_{i,j}$. These sums are all computed in the same way, and they are compulsory for checking the other computations.
\begin{gather*}
\begin{array}{lll}
  \sum_{\mathbf{r}}{D_{1,1}(\mathbf{r},\mathbf{r})} &=
    & \frac{1}{m^{2}} \times R_1^d\\
  \sum_{\mathbf{r}}{D_{2,2}(\mathbf{r},\mathbf{r})} &=
    & \left(\frac{m-1}{m}\right)^2 \times \left(\frac{m-2+_1}{(m-1)^2}\right)^d\\
  \sum_{\mathbf{r}}{D_{1,2}(\mathbf{r},\mathbf{r})} &=
    & \frac{m-1}{m^{2}}\times \left(\frac{1-R_1}{m-1}\right)^d\\
  \sum_{\mathbf{r},\mathbf{s}}{D_{1,1}} &=
    &\frac{1}{m^2} \\
  \sum_{\mathbf{r},\mathbf{s}}{D_{2,2}} &=
    &\left(\frac{m-1}{m}\right)^{2} \\
  \sum_{\mathbf{r},\mathbf{s}}{D_{1,2}} &=
    &\frac{m-1}{m^2} \\
  \sum_{\mathbf{r},\mathbf{s}}{D_{1,1}^2} &=
    &\frac{1}{m^{4}} \times  R_2^d \\
  \sum_{\mathbf{r},\mathbf{s}}{D_{2,2}^2} &=
    &\left(\frac{m-1}{m}\right)^4 \times \left(\frac{m^2-4m+6-\frac{4}{m}+R_2}{(m-1)^4}\right)^{d} \\
  \sum_{\mathbf{r},\mathbf{s}}{D_{1,2}^2} &=
    &\frac{(m-1)^2}{m^4} \times \left(\frac{1-\frac{2}{m}+R_2}{(m-1)^2}\right)^{d} \\
  \sum_{\mathbf{r},\mathbf{s}}{D_{1,1}D_{2,2}}&= &\sum_{\mathbf{r},\mathbf{s}}{D_{1,2}^2} \\
  \sum_{\mathbf{r},\mathbf{s}}{D_{1,1}D_{1,2}} &=
    &\frac{m-1}{m^4} \times \left(\frac{\frac{1}{m}-R_2}{m-1}\right)^{d} \\
  \sum_{\mathbf{r},\mathbf{s}}{D_{2,2}D_{1,2}} &=
    &\frac{(m-1)^3}{m^4} \times \left(\frac{m-3+\frac{3}{m}+R_2}{(m-1)^3}\right)^{d} 
\end{array}
,\end{gather*}
with
\begin{gather*}
 R_2=\frac{\binom{4m-3}{2m-2}}{\left((2m-1)\binom{2m-2}{m-1}\right)^2}\ \ \ 
 R_1=\frac{4^{m-1}}{(2m-1)\binom{2m-2}{m-1}}
.\end{gather*}

\end{document}